\documentclass{article}

\usepackage[margin=3cm]{geometry}
\usepackage[T1]{fontenc}
\usepackage[latin9]{inputenc}
\usepackage[english]{babel}
\usepackage[numbers]{natbib}
\usepackage{amssymb, amsfonts, verbatim}

\usepackage{mathtools}
\usepackage{graphicx}
\usepackage{pdflscape}
\usepackage{float}

\usepackage{multirow}
 
\usepackage{subfig}
\usepackage{bbm}
\usepackage{bm}
\usepackage{array}
\usepackage{bigints}
\usepackage{mathrsfs}
\usepackage[ruled]{algorithm2e}
\usepackage{lscape}

\usepackage{makecell}
\usepackage{longtable}
\usepackage{placeins}
\usepackage{textcomp}
\usepackage{hyperref}

\begin{document}
\bibliographystyle{plain}
\title{Hopping between distant basins} 

\author{Maldon Goodridge \footnote{School of Mathematics; Queen Mary, University of London; London, E1 4NS, UK. Email: ahw493@qmul.ac.uk}
\and John Moriarty \footnote{School of Mathematics, Queen Mary, University of London, London, E1 4NS, UK. Email: j.moriarty@qmul.ac.uk}
\and Jure Vogrinc \footnote{Department of Statistics, University of Warwick, Coventry, CV4 7AL, UK. Email: jure.vogrinc@warwick.ac.uk}
\and Alessandro Zocca \footnote{Department of Mathematics, Vrije Universiteit Amsterdam, 1081HV, Netherlands. Email: a.zocca@vu.nl}}

\maketitle

\begin{abstract}
We present the Basin Hopping with Skipping (BH-S) algorithm for stochastic optimisation, which replaces the perturbation step of basin hopping (BH) with a so-called skipping proposal from the rare-event sampling literature. Empirical results on benchmark optimisation surfaces demonstrate that BH-S can improve performance relative to BH by encouraging non-local exploration, that is, by hopping between distant basins. \\

\noindent \textbf{keywords}: {Basin hopping, stochastic optimisation, skipping sampler, rare events, Markov chains}
\end{abstract}

\section{Introduction and background}
\label{sec:intro}
In the literature on global optimisation of a non-convex energy landscape a source of inspiration has been methods from the theory of rare-event sampling. Examples include the methods of cross-entropy for combinatorial and continuous optimisation \cite{rubinstein1999cross} and, more recently, splitting for optimisation \cite{duan2016splitting}. In stochastic optimisation algorithms such as random search~\cite{Schumer1968}, basin hopping~\cite{leary2000global,Wales1997}, simulated annealing~\cite{Kirkpatrick1983} and the multistart method~\cite{Jain1993,Marti2003}, one or more initial points $X_0$ are perturbed in order to discover new neighbourhoods  (or `basins') of lower energy, which may then be explored by a local procedure such as gradient descent. As such algorithms discover progressively smaller energy values, the remaining lower-energy basins form a decreasing sequence of sets. Viewing the optimisation domain heuristically as a probability space and these basins as events, the discovery of smaller energy values can then also be likened to rare-event sampling. 

In this analogy, the local perturbation step plays a similar role to the proposal step in a Markov Chain Monte Carlo (MCMC) sampler (see \cite{HandbookMCMC}, \cite{RobertsRosenthal2004}). Thus in order to enhance performance, one may explore the use of alternative MCMC proposal distributions developed in the context of rare event sampling as alternative perturbation steps within stochastic optimisation routines. This is the approach we take in the present paper.

To illustrate the potential benefit of this approach consider an energy landscape having multiple, well-separated basins whose minimum energies are approximately equal to the global minimum. Then if $X_0$ lies in one such basin, separation means that local perturbations are not well suited to the direct discovery of another basin. Instead, algorithms using local perturbations to minimise over such a landscape should be non-monotonic, accepting transitions from $X_0$ to states of higher energy in the hope of later reaching lower-energy basins. In contrast, since non-local perturbation steps offer the possibility of direct moves between distant low-energy basins, they may possibly be effective on such surfaces within a monotonic optimisation algorithm. In this paper we explore the use of a particular non-local perturbation, the `skipping perturbation' of \cite{moriarty}.

Although other non-local perturbations have been proposed in the literature (see for example \cite{andricioaei2001,Pompe2020,lan2014,sminchisescu2007,sminchisescu2003mode,tjelmeland2001} in the context of MCMC), skipping has the advantage of being just as straightforward to implement as a local random walk perturbation. That is, it requires no additional information about the energy landscape beyond the ability to evaluate it pointwise.
 
We explore its use within the basin hopping (BH) algorithm~\cite{leary2000global,Wales1997}, which combines local optimisation with perturbation steps and requires only pointwise evaluations of the energy function $f$. The resulting `basin hopping with skipping' (BH-S) algorithm is thus as generally applicable as the BH algorithm.
 
The BH algorithm works as follows: the current state $X_n$ is perturbed via a random walk step to give $Y_{n}$ which is, in turn, mapped via deterministic local minimisation to a local minimum $X_{n+1}$. This local minimum point is then either accepted or rejected as the new state with probability given by the Metropolis acceptance ratio, and the procedure is repeated until a pre-determined stopping criterion is met. Due to its effectiveness and ease of implementation, the BH algorithm has been used to solve a wide array of optimisation problems (see \cite{Huang2020,Olson2012,Paleico2020} for more details).

In contrast with the non-monotonic BH algorithm, BH-S is monotonic and replaces the random walk step with a skipping perturbation over the sublevel set of the current state $X_n$. Like a flat stone skimming across water, this involves repeated perturbations in a straight line until either a point of lower energy is found, or the skipping process is halted. The BH-S algorithm, which was first outlined in \cite{moriarty}, thus provides a direct mechanism to escape local minima which contrasts with the indirect approach taken by BH. Another perspective is that BH-S alters the balance between the computational effort expended on local optimisation versus the effort spent on perturbation, typically increasing the latter while decreasing the former (cf. Table~\ref{tb:time_per_success} below).
 
Through the use of benchmark functions, the aim of the present paper is to offer guidance on tuning the method and to present a systematic overview on the types of optimisation problem on which BH-S tends to outperform BH. The rest of the paper is structured as follows: Section~\ref{sec:bhs} introduces the algorithms, empirical results are presented and discussed in Section~\ref{sec:numerics}, and Section~\ref{sec:extens} concludes.
 
\section{The BH-S algorithm}
\label{sec:bhs}
 
Consider a global optimisation problem on a rectangular subdomain $D \subset \mathbb R^d$, of the form
\begin{align}\label{eq:opt}
	\min \quad & f(\bm{x}) \quad 
	\text{ s.t. } 
	\, \bm{x} \in D: = \prod_{i=1}^d [l_i,u_i],
\end{align}
for some scalars $l_i \leq u_i$, $i=1,\dots,d$. In the rest of the paper, we will often refer to $f$ as the \textit{energy function} and to its graph as the \textit{energy landscape}. This terminology, which is similar to that of simulated annealing, is appropriate since the BH algorithm was originally conceived as a method to find the lowest energy configuration of a molecular system \cite{Wales1997}. In this section we review the BH algorithm and then introduce basin hopping with skipping (BH-S).

\subsection{Basin hopping algorithm}
\label{sub:classical}
The core idea of the basin hopping algorithm~\cite{Wales1997}, which is presented in Algorithm \ref{BHalgo}, is to supplement local deterministic optimisation by alternating it with a random perturbation step capable of escaping local minima. More specifically, inside the {\tt RandomPerturbation} procedure at step 5 of Algorithm~\ref{BHalgo} a random perturbation $W\in\mathbb R^d$ is drawn and added to the current state $X_n$ giving a state $Y_{n} = X_n+W$. Most commonly, the increment $W$ is either spherically symmetric or has independent coordinates. 
The state $Y_{n}$ becomes the starting point of a deterministic local minimisation routine. 
In our implementation of Algorithm~\ref{BHalgo} the {\tt LocalMinimisation} procedure at step 6 is performed using the limited-memory BFGS algorithm~\cite{Liu1989}, a quasi-Newton method capable of incorporating boundary constraints (although we note that other choices are possible). The resulting local minimum $U_{n}$ is then either accepted or rejected as the new state with probability equal to
\[
    \min\left(1,\exp \left( - \frac{f(U_{n}) - f(X_n)}{T} \right)\right),
\]
where $T\geq 0$ is a fixed temperature parameter. This means, in particular, that downwards steps for which $f(U_{n}) < f(X_n)$ are always accepted. The BH algorithm prescribes to repeat this basic step until a pre-defined stopping criterion is satisfied. Commonly used stopping criteria for the BH algorithm include, among others, a limit on the number of evaluations of the function $f$ or the absence of improvement over several consecutive iterations~\cite{Olson2012,Rondina2013}. The monotonic basin hopping method introduced in~\cite{leary2000global} is the BH variant corresponding to the limiting case $T=0$, in which all steps that increase the energy are rejected.

\begin{algorithm}
	\SetAlgoLined
	\everypar={\nl}
	Generate a random initial state $Y_0 \in D$\;
	$X_{0} = $ {\tt LocalMinimisation}($Y_{0}$)\;
	$n = 0$\;
	\While{Stopping criterion for $\{X_j\}_{j\leq n}$ is not satisfied:}{
		$Y_n = $ {\tt RandomPerturbation}($X_{n}$)\;
		$U_n = $ {\tt LocalMinimisation}($Y_n$)\;
		Generate $V \sim $ Uniform($[0,1]$)\;
		\uIf{V < $\min\left(1,\exp \left( - \frac{f(U_n) - f(X_n)}{T} \right)\right)$}
		    {$X_{n+1}=U_n$\;}\Else{$X_{n+1}=X_n$ \;}
        Increase $n$ by 1\;}
	
	\caption{Basin hopping}
	\label{BHalgo}
\end{algorithm}

\FloatBarrier

Basin hopping can thus be viewed as a random walk on the set of local minima of the energy landscape, which because its transition probabilities favour downhill moves to lower minima, is capable of finding the global minimum and, hence, of solving global optimisation problems. Its transition probabilities depend in a complex way on the current position, the landscape, and the perturbation step. The BH-S algorithm introduced in the next section modifies these transition probabilities, aiming to accelerate optimisation. 

\subsection{Skipping perturbations and the BH-S algorithm}
\label{sub:BHS}
In this subsection we introduce the BH-S algorithm, which differs from BH only in the perturbation step of line 5 in Algorithm \ref{BHalgo}. Instead of the random walk perturbation described above, the {\tt RandomPerturbation} procedure described in Algorithm \ref{bhmssalgo} below is applied in order to obtain $Y_{n}$. The {\tt LocalMinimisation} and acceptance steps remain identical to those in Algorithm~\ref{BHalgo}.

Given the current state $X_n$ and a fixed probability density $q$ on $\mathbb{R}^d$, the random walk perturbation of the BH algorithm can be understood as drawing a state $Y_{n}$ from the density $y \mapsto q(y-X_n)$.

In contrast the \textit{skipping perturbation} of BH-S depends on both the current state $X_n$ and a \textit{target set} $C \subseteq \mathbb{R}^d$ of states. The target set $C_n$ for the $n$-th skipping perturbation is the sublevel set of the energy function $f$ at the current point $X_n$, i.e.,
\begin{align}\label{eq:sublevelset}
    C_n := \{x \in D: f(x) \le f(X_{n})\} \subset \mathbb{R}^d.
\end{align}
    
A state $Z_1$ is drawn according to the density $q$ just as in the random walk perturbation and, if $Z_1$ does not lie in the target set $C_n$, further states $Z_2, Z_3, \ldots$ are drawn such that $X_n, Z_1, Z_2, \ldots$ lie in order on a straight line, with each distance increment $|Z_{j+1}-Z_j|$ having the same distribution as that of $|Z_1-X_n|$ conditioned on the line's direction $\frac{Z_1-X_n}{|Z_1-X_n|}$. The first state of this sequence to land in the target set $C_{n}$ becomes the state $Y_{n}$. If $C_n$ is not entered before the skipping process is halted, then $Y_{n}$ is set equal to $X_n$.

More precisely, let $x = (r,\varphi)$ be polar coordinates on $\mathbb{R}^d$ with the angular part $\varphi$ lying on the $d-1$ dimensional unit sphere $\mathbb{S}^{d-1}$. Write  $\varphi \mapsto q_{\varphi}(\varphi)$ for the marginal density of $q$ with respect to the angular part $\varphi$, which we may call the \textit{directional density} (and which we assume is strictly positive). For each $\varphi \in \mathbb{S}^{d-1}$ denote by
$$
    q_{r|\varphi}(r|\varphi):=\frac{q_{r,\varphi}(r,\varphi)}{q_{\varphi}(\varphi)}
$$
the \textit{conditional jump density}, i.e., the conditional density of the radial part $r$ given the direction $\varphi$. To construct the skipping perturbation, set $Z_0 = X_n$ and draw a random direction $\Phi \in \mathbb{S}^{d-1}$ from the directional density $q_{\varphi}$. A sequence of i.i.d.~distances $R_1, R_2, \ldots$ is then drawn from the conditional jump density $q_{r|\Phi}$, defining a sequence of modified perturbations $\{Z_k\}_{k\ge 1}$ on $\mathbb{R}^d$ by
\[
    Z_{k+1} := Z_k + \Phi R_{k+1}, \qquad k = 0,1,\ldots
\]
Since this modification of the BH perturbation is more likely to generate states $Z_k$ lying outside the optimisation domain $D$, we apply periodic boundary conditions. 

If $Z_k\in C_n$ for some $k \leq K$, where $K$ is a pre-defined maximum number of steps called the \emph{halting index}, then we set $Y_n=Z_k$ in Algorithm \ref{BHalgo} and continue to the {\tt LocalMinimisation} and acceptance steps. Alternatively if $Z_k\notin C_n$ for all $k \leq K$ we set $Y_n=X_n$. Note that although in \cite{moriarty} the halting index $K$ can be randomised, in the present setting with a known bounded domain $D$ it is sufficient to consider only fixed halting indices. 

For clarity, in the remainder of the paper we will understand the BH algorithm to mean setting $K=1$ in Algorithm~\ref{bhmssalgo}. In all simulations we set the perturbation $q$ to be a spherically symmetric and Gaussian with standard deviation $\sigma$, although other choices are possible (see the discussion in Section~4.1 of \cite{moriarty}). In the next section we explore for which types of energy function $f$ BH-S offers an advantage over BH, and also discuss the choice of halting index.

\begin{algorithm}[!h]
\SetKwInOut{Input}{Input}\SetKwInOut{Output}{Output}
    \SetAlgoLined
	\Input{State $X_n \in \mathbb{R}^d$}
    \Output{Randomly perturbed state $Y_n \in \mathbb{R}^d$}
	\everypar={\nl}
Set $Z_0=X_n$ \;
Generate an initial perturbation $W$ distributed according to the density $w \mapsto q(w-X_n)$ \;
Calculate the direction 
$$\Phi \quad=\quad \frac{(W - X_n)}{\|W - X_n\|}\,;$$
Set $k=1$ and $Z_1:=W$\;
\While{$f(Z_k) > f(X_n)$ \textbf{\textup{and}} $k< K$}{
	Generate an independent distance increment $R$ distributed as $\|W-X_n\|$ given $\Phi$ \;
	Set $Z_{k+1}=Z_k+ \Phi R$ \;
	Increase $k$ by one \;
}
Set $Y_n:=Z_k$\;
\caption{{\tt RandomPerturbation} subroutine for BH-S} 
\label{bhmssalgo}
\end{algorithm}

\FloatBarrier

\section{Empirical results}
\label{sec:numerics}
In this section we aim to explore on which types of optimisation problem BH-S tends to outperform BH and {\it vice versa} using a set of benchmark energy landscapes with known global minima from \cite{optonline,optlit,optsfu}. To facilitate discussion of landscape geometry we initially restrict attention to two-dimensional energy functions, before considering higher dimensions in Section \ref{sub:dimen}.

In Subsection~\ref{sub:distbas} we show that, if an energy landscape has \textit{distant basins} (recall that with the word `basin' we refer to the neighbourhood of a local minimum) then BH-S tends to offer an advantage. Otherwise, as described in Subsection~\ref{sub:nondistbas}, BH is to be preferred since any benefit from BH-S is then typically outweighed by its additional computational overhead. We also explore the effect of the state space dimension $d$ on the performance of both algorithms and offer guidance on tuning the BH-S method, including strategies to improve exploration of challenging energy landscapes.

\subsection{Methodology}
\label{sub:method}
For each benchmark energy landscape, we compare the performance of BH-S to that of BH with temperature $T = 1$, in both cases taking the density $q$ of the initial perturbation as the centred Gaussian
\begin{align}\label{eq:propq}
q \sim \mathcal{N}(\bm{0},\sigma^2 \cdot I_d),
\end{align}
where $I_d$ is the $d \times d$ identity matrix and the parameter $\sigma$ allows for tuning, as follows. Both the BH and BH-S algorithms are run on a set of uniformly distributed initial states $I:=\{X_0^{(n)} \in D, \, n=1,\ldots,|I|\}$. These initial states are used sequentially until the computational budget of 300 seconds of CPU time has elapsed, and the corresponding set of final states is recorded. To account for numerical tolerance, we consider a run to to have successfully identified the global minimum $x^* \in D$ if its final state lies in $\mathcal{G}:=\{x \in \mathbb{R}^d~:~||x-x^*||\le 10^{-5}\}$ (this choice excludes all non-global minima for all benchmark landscapes). 

The performance of each algorithm is then assessed with respect to two metrics:
\begin{itemize}
    \item {\bf Reliability,} defined as the proportion of runs terminating in $\mathcal{G}$, 
    \item {\bf Efficiency,} defined as the number of runs terminating in $\mathcal{G}$.
\end{itemize}
 
We write $\rho_c$ and $\rho_s$ for the reliability of the BH and BH-S algorithms respectively, while $\epsilon_c$ and $\epsilon_s$ denote their respective efficiencies. The BH and BH-S algorithms are individually tuned for each function by selecting $\sigma$ and $K$ to maximise their efficiency.

In order to understand the role played by the skipping perturbation, we also record diagnostics on the average size of perturbations. For each new state $X_{n+1} \neq X_n$ accepted in Algorithm \ref{BHalgo}, define the \textit{perturbation distance} $J$ as $\|Y_n-X_n\|$, the Euclidean distance between the state $X_n$ at step $n$ and its perturbation $Y_n$. For each run of an algorithm, the mean $\overline{J}$ of these perturbation distances is recorded. Then for each 300 second budget, the \textit{expected mean jump distance} $\upsilon$ is the average $\upsilon:=N^{-1} \sum_{n=1}^N \overline{J}^{(n)}$, where $N$ is the number of runs realised within the time budget. For the BH-S algorithm, $\upsilon$ is calculated separately for the accepted random walk perturbations (that is, those for which $Y_n=Z_1$ in Algorithm \ref{bhmssalgo}) and the accepted skipping perturbations (those for which $Y_n=Z_k$ with $k \geq 2 $ in Algorithm \ref{bhmssalgo}), denoting these by $\upsilon_1 \text{ and } \upsilon_s$ respectively. 

The simulations were conducted on a single core using Python 3.7, using the {\tt basinhopping} routine in SciPy version 1.6.2 for the BH algorithm. Results for all considered landscapes are presented in the Appendix. Jupyter notebook files used to conduct simulations are available at \url{https://github.com/ahw493/Basin-Hopping-with-Skipping}.

\subsection{Exploratory analysis}
\label{sub:opt_surface}

\begin{figure}[h]
 
	\hspace{-1cm}
	\includegraphics[width=1.15\textwidth]{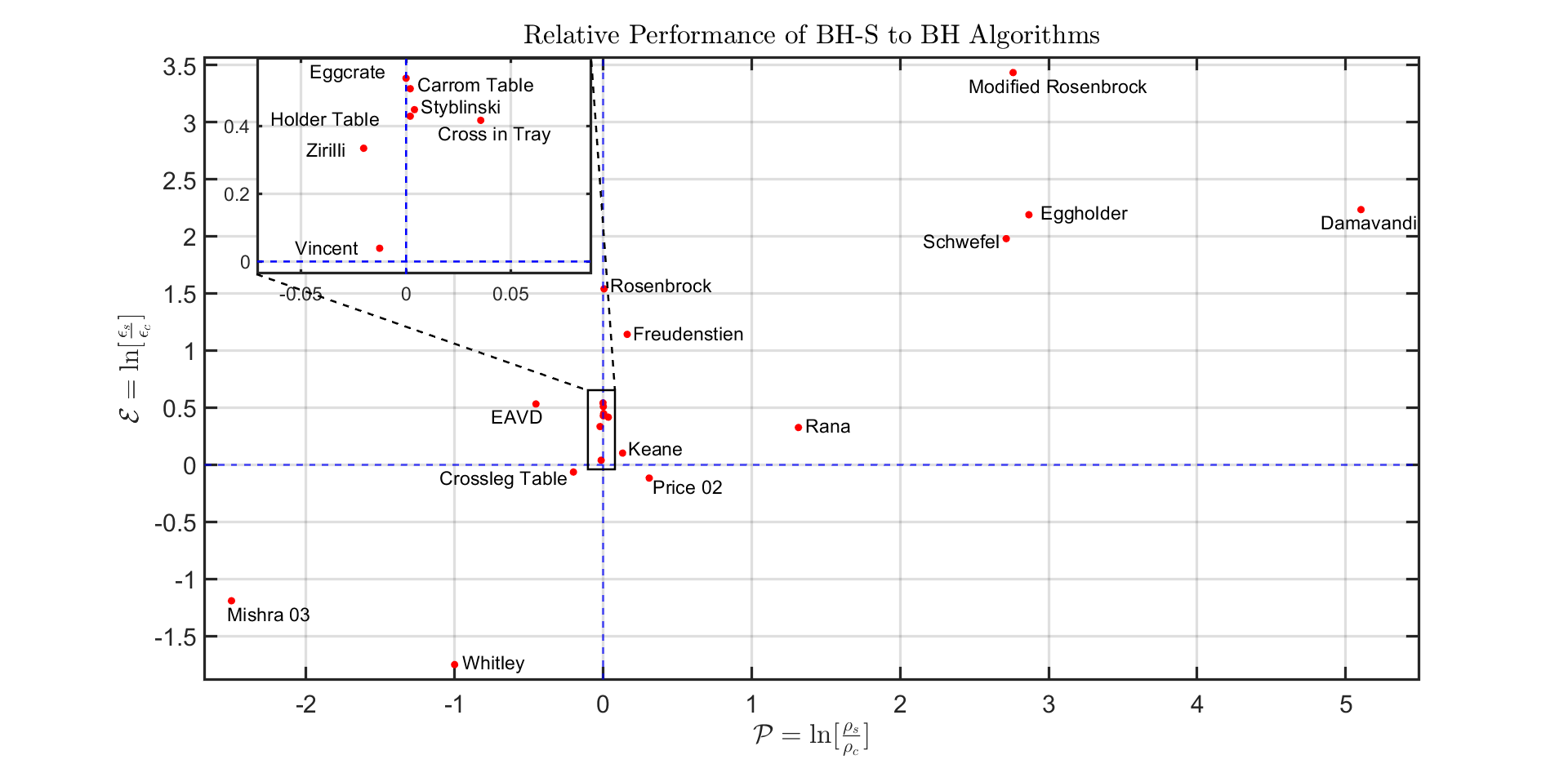}
	\caption{Scatterplot of relative efficiency $\mathcal{E}=\ln(\epsilon_s /\epsilon_c)$ versus relative reliability $\mathcal{P}=\ln(\rho_s /\rho_c)$ for the BH and BH-S algorithms on benchmark energy landscapes}
	\label{fg:pibeps}
\end{figure}

As an exploratory comparison between BH and BH-S, their relative efficiency $\rho_s /\rho_c$ and reliability $\epsilon_s /\epsilon_c$ were calculated for each benchmark energy landscape and plotted in Figure~\ref{fg:pibeps}. 

Landscapes in the first quadrant of Figure~\ref{fg:pibeps} represent cases where the BH-S algorithm exhibits both greater reliability and greater efficiency than BH. The common feature among these landscapes, which are plotted in the Appendix, might be called {\it distant basins}: that is, basins separated by sufficient Euclidean distance that the random walk performed by the BH algorithm is unlikely to transition directly between them. While indirect transitions between such basins may be  possible, they require a suitable combination of steps to be made. Such indirect transitions may carry significant computational expense, for example if suitable combinations of steps are long or relatively unlikely. In the BH-S algorithm, by contrast, the linear sequence of steps taken by the skipping perturbation enables direct transitions even between distant basins.

Conversely, landscapes lying in the lower-left quadrant of Figure~\ref{fg:pibeps} represent cases where the BH-S algorithm is both less reliable and less efficient than BH. As explored more extensively later in Subsection \ref{sub:nondistbas}, for each of these landscapes, if the energy of the state $X_n$ is close to the global minimum value $f(x^*)$ then the corresponding sublevel set $C_n$ has almost zero volume. This means that even if the skipping perturbation  traverses the distance between basins, the states $Z_1,\ldots,Z_k$ are unlikely to fall in $C_n$ due to its small volume. Since the BH algorithm is non-monotonic, it does not suffer from the same issue and outperforms BH-S for these landscapes.

Figure~\ref{fg:pibeps} displays a positive correlation between relative efficiency and relative reliability. However for several landscapes (which lie near the vertical axis) the performance of BH and BH-S cannot be clearly distinguished on the basis of reliability alone. As confirmed by the Appendix, this is typically because both algorithms have reliability close to 100\%. Nevertheless the algorithms differ in their efficiency, with BH-S observed to be more efficient than BH for each such landscape. One surface also lies in each of the second and fourth quadrants.

\setlength\tabcolsep{2pt}
\renewcommand{\arraystretch}{1.5}
\begin{table}[!ht]
    \centering
	\caption{CPU time spent on the perturbation and local minimisation steps by the BH and BH-S algorithms for the test functions discussed in Sections~\ref{sub:distbas}--\ref{sec:sc}, normalised by efficiency.}

	\begin{tabular}[h!]{llcccc}		
		\hline\noalign{\smallskip}
        &&\multicolumn{4}{c}{$\frac{\text{Time spent}}{\text{Efficiency}}$, s }\\
		&&\multicolumn{2}{c}{\makecell{BH}}& \multicolumn{2}{c}{\makecell{BH-S}} \\ \noalign{\smallskip}\hline\noalign{\smallskip}
		
		Location in Figure~\ref{fg:pibeps}&Function & \makecell{Perturbation} & \makecell{Local \\Minimisation} & \makecell{Skipping \\Perturbation} & \makecell{ Local\\ Minimisation}\\ 
		
		\noalign{\smallskip}\hline\noalign{\smallskip}
		
	\multirow{2}{*}{\makecell{First quadrant \\(Section~\ref{sub:distbas})}}&	Egg-holder 	 & 	 0.72	 & 	6.76 	& 	0.73  & 	 0.19 \\ \cline{2-6}
	& Modified Rosenbrock 	 & 	 0.46	 & 	13.53 	& 	0.22  & 	 0.10 \\ \hline
     \multirow{2}{*}{\makecell{Third quadrant\\ (Section~\ref{sub:nondistbas})}} & Mishra-03 	 & 	 0.02	 & 	1.75 	& 	1.74  & 	 3.21 \\ \cline{2-6}
	&Whitley 	 & 	 0.01	 & 	0.68 	& 	4.92  & 	 1.58 \\ \hline
    \multirow{2}{*}{\makecell{Special cases\\ (Section~\ref{sec:sc})}} &   Rosenbrock	 & 	 0.01	 & 	0.13 	& 	0.05  & 	 0.01 \\ \cline{2-6}
      &  Styblinksi	 & 	 0.02	 & 	0.06 	& 	0.06  & 	 0.01 \\

		\noalign{\smallskip}\hline
	\end{tabular}
	\label{tb:time_per_success}
\end{table}
 
Further exploratory analysis is provided in Table~\ref{tb:time_per_success}, which indicates average CPU time spent on the perturbation versus the local minimisation steps for each algorithm. To facilitate comparisons between the two algorithms, in each case the total time spent is normalised by the algorithm's efficiency (as defined in Section~\ref{sub:method}). This demonstrates that the BH algorithm invests a large majority of processor time in the local minimisation step, with relatively little devoted to the perturbation step. While the ratio between processor time spent on local minimisation and perturbation is more problem-dependent for BH-S, the balance appears to be shifted in favour of perturbation.

The BH-S perturbation step is more expensive by construction, since it requires between $1$ and $K$ evaluations of the energy function $f$ (depending on the sublevel set of the current state), whereas each BH perturbation requires just one evaluation of $f$. However in Table~\ref{tb:time_per_success}, for the Damavandi, Schwefel, Modified Rosenbrock and Egg-holder functions for which BH-S works well (cf. Figure~\ref{fg:pibeps}), after normalisation the BH-S algorithm spent approximately the same or less CPU time than BH on perturbation, in addition to spending less time on local minimisation. Thus for these landscapes which favour BH-S, perturbation steps were not only less frequent for BH-S (again, after normalisation by efficiency) than BH, but the monotonic BH-S perturbations also reduced the total computational burden arising from the local minimisation step. 

Conversely it was noted above that for landscapes in the third quadrant of Figure~\ref{fg:pibeps}, if the energy of the state $X_n$ is close to the global minimum value $f(x^*)$ then the corresponding sublevel set $C_n$ has almost zero volume. This represents the worst case for the BH-S perturbation: if the states $Z_1,\ldots,Z_k$ all lie outside the sublevel set then the perturbation requires the maximum number $k$ of evaluations of the energy function, but nevertheless the perturbed state $Y_n$ is rejected and $X_{n+1}=X_n$, so the optimisation procedure does not advance. Indeed for the Mishra-03 and Whitley functions in Table~\ref{tb:time_per_success}, the efficiency normalised CPU time invested in perturbations is two orders of magnitude greater for BH-S than for BH. For these landscapes, the efficiency normalised computational burden from local minimisation is also observed to be greater for BH-S than for BH, although the reasons for this are less clear.

Guided by the exploratory analysis of Figure \ref{fg:distbas}, in Sections \ref{sub:distbas}--\ref{sec:sc} we study both algorithms' performance on specific energy landscapes in greater detail.

\subsection{Landscapes favouring the BH-S algorithm}
\label{sub:distbas}

\begin{figure}[h!]

	\centering
	\subfloat[][Modified Rosenbrock function]{\includegraphics[width=0.45\textwidth]{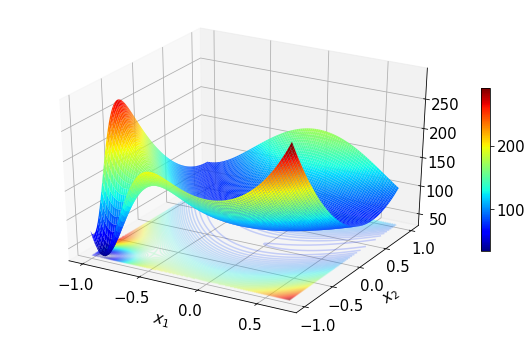}\label{fg:modros-graph}}\qquad
	\subfloat[][Egg-Holder function]{\includegraphics[width=0.45\textwidth]{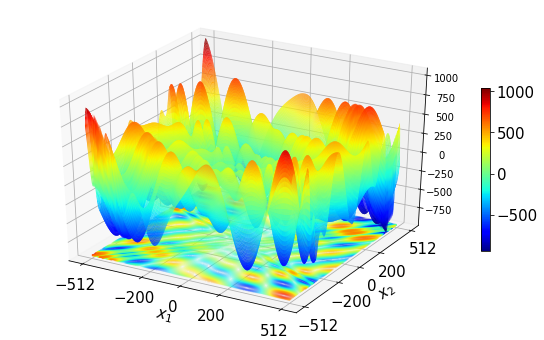}\label{fg:egg-graph}} \\
	\vspace{3mm}
	\subfloat[][\makecell{Sublevel set $C =\{x\in \mathbb R^2 ~:~ f(x)< 100\}$\\of the Modified Rosenbrock function}]{\includegraphics[width=0.45\textwidth]{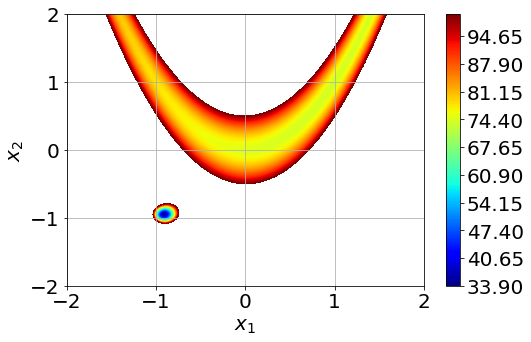}\label{fg:modros-graphsub}}\quad
	\subfloat[][ \makecell{Sublevel set $C = \{x \in \mathbb R^2 ~:~ f(x) < -700\}$\\ of the Egg-Holder function}]{\includegraphics[width=0.45\textwidth]{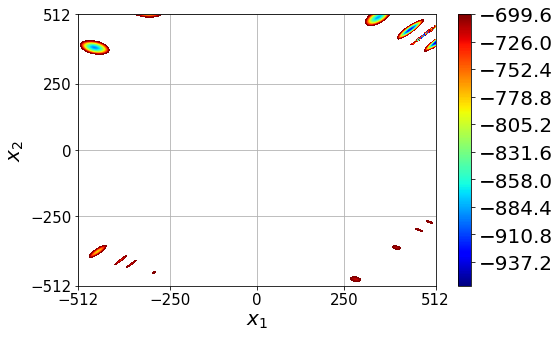}\label{fg:egg-graphsub}}  
	\caption{Examples of energy landscapes from the first quadrant of Figure~\ref{fg:pibeps} }
	\label{fg:distbas}
\end{figure}

Figure \ref{fg:distbas} plots two landscapes from the first quadrant of Figure~\ref{fg:pibeps}--that is, landscapes which favour the BH-S algorithm over BH. For each landscape,  the sublevel set of a level above the global minimum $f(x^*)$ is also plotted. The Modified Rosenbrock energy function is given by

\[
    f(x) =  74 +100(x_2 - x_1^2)^2 + (1 - x_1)^2 - 400 \exp\bigg[-\frac{(x_1 +1)^2 +(x_2 +1)^2}{0.1}\bigg],
\]
and we take the domain $D = [-2,\ 2]^2$, with global minimum $x^* = (-0.95, \ -0.95)$ \cite{optonline}.

The Egg-Holder energy function is 

\[
    f(x) = -(x_2 +47)\sin{\bigg( \sqrt{|x_2 + \frac{x_1}{2}+47|}  \bigg) }-x_1 \sin{\bigg( \sqrt{ |x_1 - (x_2+47)|} \bigg)},
\]
and we take the domain $D = [-512,512]^2$, with global minimum at $x^{*} = (512, 404.2319)$ \cite{optlit}.

From Figure~\ref{fg:modros-graph}, the modified Rosenbrock function has two basins: a larger basin with a U-shaped valley and a smaller, well-shaped basin. To transition from the minimum of the valley to the minimum of the well, the BH algorithm would require a relatively large perturbation step landing directly in the well, otherwise the local optimisation procedure would take it back to the minimum of the valley. Even for an optimal choice of $\sigma$, which would require {\it a priori} knowledge about the landscape, such perturbations would be unlikely.

In contrast, if the initial point $X_0$ lies at the minimum of the valley, the BH-S algorithm aims to skip across the domain and enter its sublevel set $C_0$ as defined in \eqref{eq:sublevelset}. From Figure~\ref{fg:modros-graphsub}, this will correspond to entering an approximately  circular basin near the point $(-1,-1)$ in the domain. By Algorithm \ref{bhmssalgo}, the skipping perturbation has the potential to enter that basin provided that the straight line issuing from $X_0$ in the initial direction $\Phi$ in Algorithm \ref{bhmssalgo} intersects it. In particular, this ability is robust to the choice of standard deviation $\sigma$ provided that the halting index $K$ is chosen appropriately (see the discussion on tuning in Subsection~\ref{sub:tuning} below). 
 
\begin{table} 
\centering
\caption{Reliability $\rho$ and efficiency $\epsilon$ of the BH and BH-S algorithms, with and without the application of periodic boundary conditions, for the test functions discussed in Sections~\ref{sub:distbas}--\ref{sec:sc}.}
\begin{tabular}[h!]{cllcccccc}	

		\hline\noalign{\smallskip}
                 
		Location in Figure~\ref{fg:pibeps}&Function&\makecell{Boundary Condition}&$\rho_c$& $\epsilon_c$&\hspace{0.5cm}&$\rho_s$&$\epsilon_s$\\ 
		\noalign{\smallskip}\hline\noalign{\smallskip}
		
\multirow{4}{*}{\makecell{First quadrant\\ (Section~\ref{sub:distbas})}}&\multirow{2}{*}{Egg-Holder}	&	Periodic    	&	1.8\%	    &	77	&\hspace{0.5cm}	&	29.4\%	&	280	\\
                        	&&	Non-periodic	&	1.4\%	    &	58	&	&	81.7\%	&	815	\\\cline{2-8}
&\multirow{2}{*}{Modified Rosenbrock}         	&	Periodic   	&	18.8\%	&	783 	&	&	100\%	&	1249	\\
&                        	&	Non-periodic	&	17.2\%	&	597 	&	&	96.7\%	&	321	\\\hline
\multirow{4}{*}{\makecell{Third quadrant\\ (Section~\ref{sub:nondistbas})}}&\multirow{2}{*}{Mishra-03}	&	Periodic    	&	77.9\%	&	811	    &	&	6.2\%	&	78	\\
&                        	&	Non-periodic	&	49.6\%	&	474	    &	&	68.4\%	&	576	\\\cline{2-8}
&\multirow{2}{*}{Whitley}	&	Periodic    	&	90.1\%	&	137	    &	&	13.9\%	&	33	\\
&                        	&	Non-periodic	&	86.8\%	&	138	    &	&	13.4\%	&	34 \\ \hline
\multirow{4}{*}{\makecell{Special cases\\ (Section~\ref{sec:sc})}}&\multirow{2}{*}{Rosenbrock}	&	Periodic    	&	100\%	    &	942	&	&	100\%	&	1376	\\
&                        	&	Non-periodic	&	100\%	    &	1006	&	&	100\%	&	1524	\\\cline{2-8}
&\multirow{2}{*}{Styblinski}	&	Periodic    	&	99.7\%   	&	366	&	&	99.9\%	&	932	\\
&                        	&	Non-periodic	&	99.4\%	&	526	    &	&	92.3\%	&	784	\\	
\noalign{\smallskip}\noalign{\smallskip}\hline
\noalign{\smallskip}\hline
	\end{tabular}
	\label{tb:egg_boundaries}
\end{table}
\FloatBarrier

From Figure~\ref{fg:egg-graph} the Egg-Holder function has multiple basins, many of which have near-global minima. Figure~\ref{fg:egg-graphsub} shows that the deepest basins lie in four groups, one group per corner of the domain. Within each group, the basins are close in the Euclidean distance and so perturbations are likely to enter different basins within that group. Also, the basins in each group have similar depths (that is, similar local minimum energies), making the acceptance ratio in Algorithm \ref{BHalgo} high for such within-group perturbations. As a result the BH algorithm is likely to walk regularly between within-group local minima.
Also from Figure~\ref{fg:egg-graph}, the Egg-Holder function has shallower basins distributed throughout its domain. As discussed in Subsection \ref{sub:opt_surface} these provide an indirect, although potentially computationally expensive, route for BH to cross between the four groups of Figure~\ref{fg:egg-graphsub}. 

However between groups the Euclidean distance is large, creating the same challenge for BH as with the modified Rosenbrock function: even for optimally chosen $\sigma$, which would require {\it a priori} knowledge of the landscape, transitions between groups are relatively rare.

In contrast, the BH-S algorithm is capable of moving between the four groups in Figure~\ref{fg:egg-graphsub} provided the initial direction $\Phi$ of its skipping perturbation intersects a different group. The likelihood of such an intersection is increased by both the length of the skipping chain and the use of periodic boundary conditions in the BH-S algorithm, and is again robust with respect to the choice of standard deviation $\sigma$. 

Regarding the application of periodic boundary conditions to the domain $D$, we have argued that they are natural for BH-S, since otherwise long skipping chains would tend to exit the domain $D$. In contrast, they are not implemented for the BH algorithm in the results of Figure~\ref{fg:pibeps} and Table~\ref{tb:results}. One may therefore ask whether it is their use, rather than the skipping perturbation of BH-S, which yields any observed improvement. To explore this, Table \ref{tb:egg_boundaries} illustrates the effect of imposing periodic boundary conditions on the performance of both the BH and BH-S algorithms. Interestingly the performance of {BH-S} on the Egg-Holder landscape is improved without their use (a fact which appears to be driven by the proximity of its global minimiser $x^*$ to the boundary). In general, it is clear from Table~\ref{tb:egg_boundaries} that for both algorithms their benefit or disbenefit is problem-dependent and the skipping perturbation explains a distinct and material part of the observed improvements relative to BH.

It can be observed from Table~\ref{tb:results} that the expected mean jump distances $\upsilon_s$ and $\upsilon_c$ (defined in Subsection \ref{sub:method}) typically satisfy $\upsilon_s >> \upsilon_c$ for landscapes in the first quadrant of Figure~\ref{fg:pibeps}. This confirms quantitatively the success of BH-S in hopping between distant basins. The cost of this feature is that the BH-S skipping perturbation is more computationally intensive than the random walk perturbation of BH.

Without skipping (that is, using the halting index $K=1$ in Algorithm~\ref{bhmssalgo}), BH-S would reduce to the monotonic basin hopping method of~\cite{leary2000global} and the {\em initial} perturbation $W$ of Algorithm \ref{bhmssalgo} would simply be either accepted or rejected. One may therefore also ask whether this increase in the expected mean jump distance is induced by the skipping mechanism of BH-S, or simply by its monotonicity. To address this, recall that Algorithm~\ref{bhmssalgo} first perturbs the current state $X_n$ to give an initial perturbation $Z_1:=W$. Then if $f(W)>f(X_n)$, the initial perturbation is modified to $Z_2$, and so on, until either a state $Z_k$ is generated with $f(Z_k)\leq f(X_n)$ or skipping is halted. If such a $Z_k$ is found then it may be accepted by setting $X_{n+1}=Z_k$ or rejected. The Appendix records the proportion of accepted BH-S perturbations $X_{n+1} = Z_k$ for which $k>1$. Indeed, for many landscapes in the first quadrant of Figure~\ref{fg:pibeps} this proportion is 100\%. That is, for such landscapes, each accepted perturbation $X_{n+1}$ required the skipping mechanism since none of the initial perturbations had lower energy than the current state $X_n$.

\subsection{Landscapes favouring the BH algorithm}
\label{sub:nondistbas}
\begin{figure}[!ht]
	\hspace*{-1.5cm}
	\centering
	 \qquad
	 
	 \subfloat[][Mishra-03 Function]{\includegraphics[width=0.45\textwidth]{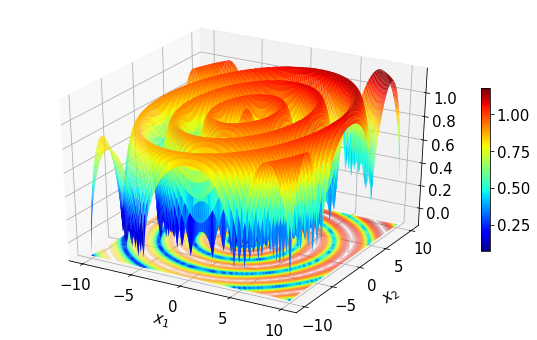}\label{fg:mish_graph} } \quad
	  \subfloat[][Whitley Function]{\includegraphics[width=0.45\textwidth]{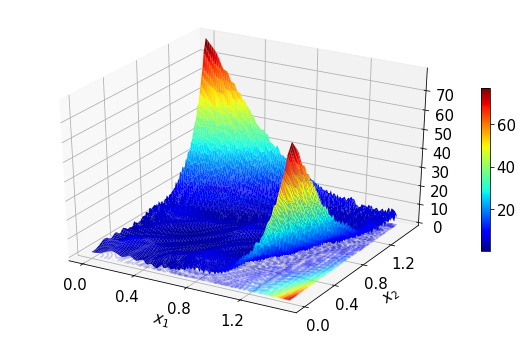}\label{fg:whit_graph}} \\
     \par\bigskip
      \subfloat[][Sublevel set $C=\{x:f(x)<1\}$ of the Mishra-03 function]{\includegraphics[width=0.45\textwidth]{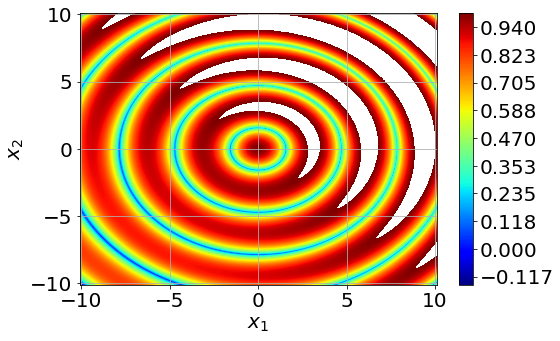}\label{fg:mish_L1}}\quad
	  \subfloat[][Sublevel set $C=\{x:f(x)<5\}$ of the Whitley function]{\includegraphics[width=0.42\textwidth]{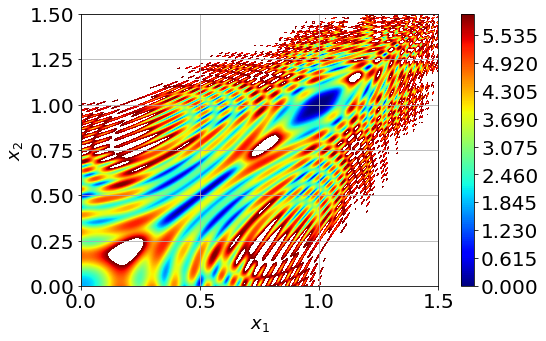}\label{fg:whit_L5}}  \\
     \par\bigskip
    \subfloat[][Sublevel set $C=\{x:f(x)<-0.01\}$ of the Mishra-03 function]{\includegraphics[width=0.45\textwidth]{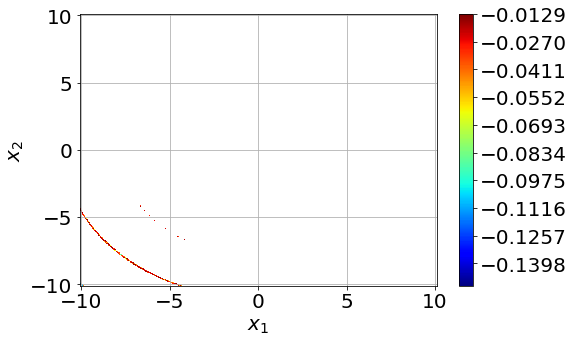}\label{fg:mish_L01}} \quad
    \subfloat[][Sublevel set $C=\{x:f(x)<0.1\}$ of the Whitley function]{\includegraphics[width=0.45\textwidth]{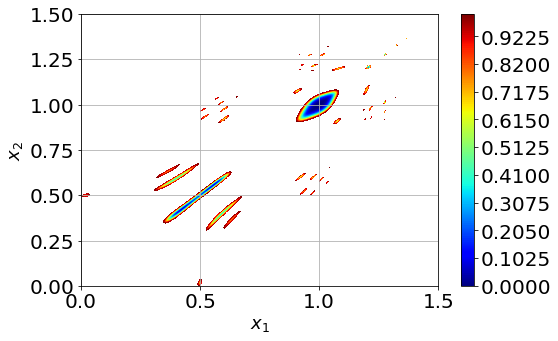}\label{fg:whit_L02}}	
 	 
	\caption{Examples of energy landscapes favouring the BH algorithm} 
	\label{fg:nondistbas}
\end{figure}

Figure~\ref{fg:nondistbas} plots two landscapes from the third quadrant of Figure~\ref{fg:pibeps}, on which the BH algorithm outperforms BH-S, each with two sublevel sets above the global minimum $f(x^*)$. The Mishra-03 function is 

\[
    f(x) :=  \sqrt{|\cos\big(\sqrt{|x_1^2 + x_2^2|}\big)|} + 0.01(x_1 + x_2),
\]
and the domain $D = [-10,10]^2$ gives $x^* = (-8.466,-10)$  \cite{optlit}. The Whitley function $f: \mathbb{R}^2 \to \mathbb{R}$, given by

\[
    f(x) := \sum_{i=1}^2 \sum_{j=1}^2 \bigg ( \frac{\big(100(x_{i}^2 - x_j)^2 + (1 - x_j)^2\big)^2}{4000} - \cos\big(100(x_{i}^2 - x_j)^2 + (1 -x_j)^2 +1 \big)\bigg ),
\]
has global minimum $x^* =(1,1)$ on the domain $[0,1.5]^2$ \cite{optlit}.

From Figure~\ref{fg:mish_graph}, the Mishra-03 function is highly irregular and has many basins which appear almost point-like. Figure~\ref{fg:mish_L01} confirms that the situation outlined in Section~\ref{sub:opt_surface} applies to this landscape. That is, for states $X_n$ with energy close to the global minimum value $f(x^*)$, the corresponding sublevel set $C_n$ has almost zero volume and the states $Z_1,Z_2,\ldots,$ of Algorithm~\ref{bhmssalgo} are unlikely to fall in $C_n$. 

The deepest basins of Mishra-03 form groups arranged in concentric circular arcs. Since the Euclidean distances both within and between these groups are relatively small, the BH algorithm is able to move frequently both within and between groups without requiring precise tuning of the standard deviation parameter $\sigma$. In particular, it outperforms BH-S on this landscape. 

Similarly from Figure~\ref{fg:whit_L5}, the deepest basins of the Whitley function can be seen either as forming one group, or as a small number of groups close to each other. Thus, as for Mishra-03, the BH algorithm is able to move frequently between them while nevertheless being robust to the choice of the standard deviation parameter $\sigma$. As with Mishra-03, however, from Figure~\ref{fg:whit_L02} the sublevel sets $C_n$ corresponding to near-global minimum states $X_n$ have low volume. Thus it is more challenging for BH-S to transition between the deepest basins, and BH outperforms BH-S on this landscape.


These limitations of the BH-S routine can be mitigated by alternating between a monotonic and non-monotonic perturbation step. In Subsection~\ref{sub:alt} we provide a discussion on how this alternating perturbation can be implemented.

\subsection{Special cases}\label{sec:sc}

It was noted in Section \ref{sub:opt_surface} that for several landscapes lying near the vertical axis, both BH and BH-S algorithms have reliability close to 100\%. For these surfaces BH-S typically has greater efficiency simply because of its monotonicity, since no further computational effort is expended on local optimisation once the global optimum is reached. The Holder Table and Carrom Table landscapes have multiple distant `legs', each leg being the basin of a global minimum point. In this case, the ability of BH-S to skip between distant basins is not reflected in either its efficiency or its reliability, although it would clearly be beneficial if the goal was to identify the number of global minima in the landscape. 

\subsection{Scaling with dimension}
\label{sub:dimen}
 
In this section we aim to illustrate the performance of the BH-S algorithm as the dimension of the optimisation problem increases. For this we focus on  Schwefel-07, a landscape with `distant basins' which is also defined for higher dimensions. It is given by the function $f_d: \mathbb{R}^d \to \mathbb{R}$, where
\[
    f_d(x) = 418.9829\times d - \sum_{i=1}^d x_i \sin(\sqrt{|x_i|}),
\]
and has global minimum $(421.0)^d$ on the domain $D_d=[-500,500]^d$ 
~\cite{optlit}. 

\begin{figure}[!ht]
\centering
    \subfloat[][Percent of samples in desired basin:  $\rho$]{\includegraphics[width=0.45\textwidth]{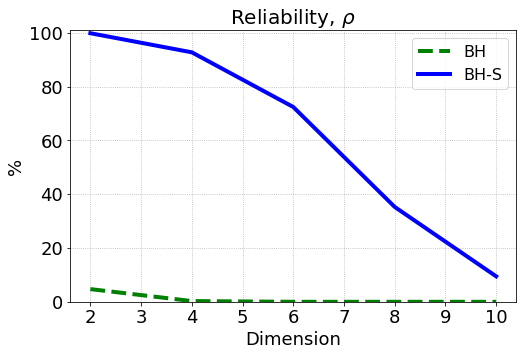}\label{fg:sch_per_inputA}} \qquad
    \subfloat[][Total trials in desired basin: $\epsilon$ ]{\includegraphics[width=0.45\textwidth]{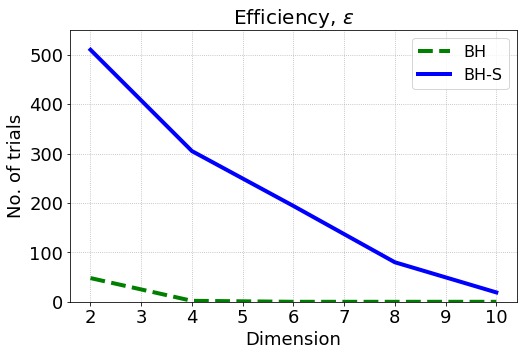}\label{fg:sch_per_inputB}}\\

	\caption{Comparison of BH and BH-S performance when applied to the Schwefel-07 function while varying the dimension $d$ of the domain $D$. We set $\sigma = 20$ for both algorithms and the BH-S has halting index $K = 50$. These parameters were close to optimal for both algorithms. Each simulation used a CPU time budget of 300s.}
	\label{fg:sch_per_input}
\end{figure}

From Figure~\ref{fg:sch_per_inputA}, the reliability and also the efficiency of both algorithms decrease approximately linearly with increased dimension. Recall that relative to BH, the strength of BH-S lies in its ability to transition directly between distant basins. From Algorithm~\ref{bhmssalgo}, in order to transition directly to the global minimum basin, it is necessary for the line from the current state in the random direction $\Phi$ to intersect that basin. As $\Phi$ is drawn from a space of  dimension $d-1$, heuristically this becomes less likely as $d$ increases.

In contrast, the BH algorithm should rely to a greater extent on indirect transitions from its current state to the global minimum. By statistical independence, the probability of a particular indirect transition is the product of the probabilities of its constituent steps. Since the probability of each step decays with dimension as discussed above for BH-S, this suggests that the performance of BH will degrade more rapidly with dimension than BH-S.

This is borne out in Figure~\ref{fg:sch_per_inputA}, where BH fails to locate $x^{*}$ within the 300 second budget for any dimension $d\ge4$, while BH-S continues to locate $x^*$ (albeit with decreasing reliability and efficiency) until dimension $d=11$. Indeed, the reliability of BH-S for this landscape is above 50\% for dimensions $d \le 7$. 

\subsection{Tuning}
\label{sub:tuning}

Both BH and BH-S have the parameter $\sigma$, the standard deviation of the centred Gaussian density $q$ used to generate the initial perturbation. As noted above, the initial perturbation is analogous  to a Metropolis-Hastings (MH) proposal in MCMC. The MH literature highlights the importance of tuning such proposals, guided either by theory or by careful experimentation~\cite{Gamerman2006,Metropolis1953}. Following this analogy, in this section we explore the choice of $\sigma$ and also of the BH-S halting index $K$. To facilitate this discussion we restrict attention to the two-dimensional Egg-Holder function.

\begin{figure}[!h]
    \centering
    \subfloat[][Percentage of trials which successfully reported the correct global minimum]{\includegraphics[width=0.45\textwidth]{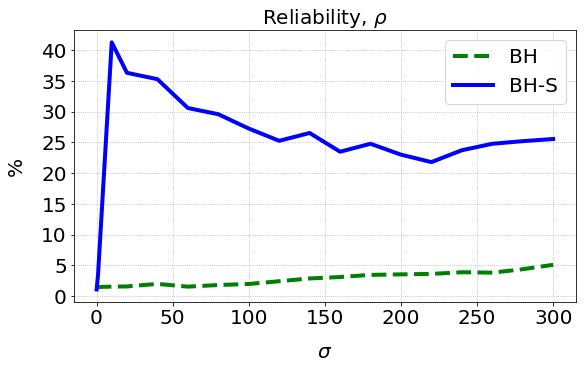} \label{fig:eggresultsA}} \quad
    \subfloat[][Total number of trials which successfully reported the correct global minimum.]{\includegraphics[width=0.46\textwidth]{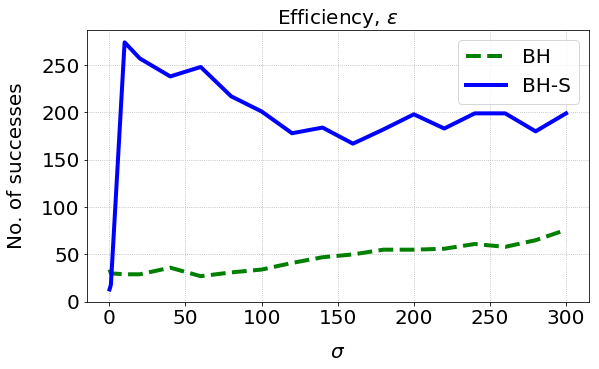}\label{fig:eggresultsB}}\\
    
    \subfloat[][Total perturbation steps conducted  during the 300s time budget.]{\includegraphics[width=0.45\textwidth]{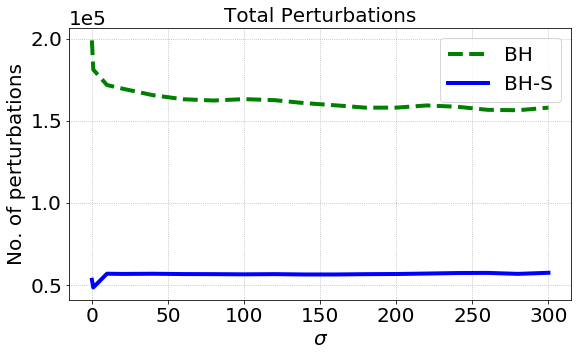}\label{fig:eggresultsC}} 

    \caption{Comparison of individually tuned BH-S and BH performances on the Egg-holder function. Set-up: CPU time budget of 300 seconds; stopping criteria: 50 perturbations; the halting index for skipping perturbation is set to $K = 25$ for all simulations.} 
    \label{fg:eggresults}
\end{figure}

Figure~\ref{fg:eggresults} plots the reliability and efficiency of both BH and BH-S as $\sigma$ varies between 0 and 300 (recall that the domain $D=[-512,512]^2$; also, we set $K=25$ for BH-S). Clearly, for both algorithms $\sigma$ should not be very small ($\leq 10$). In that case the random walk step $W$ is likely to land in the same basin as the current point $X_n$, so that the local optimisation step maps the perturbation back to $X_n$ and the algorithms do not advance.

We note first from Figure~\ref{fg:eggresults} that both the reliability and efficiency of the BH algorithm increase approximately linearly within this range as $\sigma$ increases. As discussed in Section~\ref{sub:distbas}, this reflects the fact that as $\sigma$ increases, direct transitions between the four groups of deepest basins become more likely. In contrast, and again confirming the discussion in Section~\ref{sub:distbas}, both the efficiency and reliability of BH-S appear to be rather robust to the choice of $\sigma$. 

\begin{figure}[!ht]
	\centering
    \subfloat[][Percent of samples in desired basin:  $\rho_s$]{\includegraphics[width=0.45\textwidth]{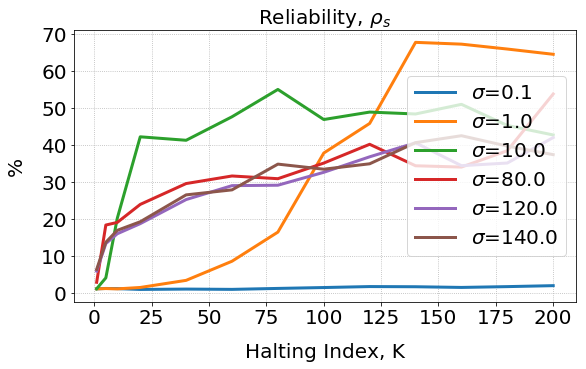}\label{fg:egg_ksenA}} \qquad
    \subfloat[][Total trials in desired basin: $\epsilon_s$ ]{\includegraphics[width=0.45\textwidth]{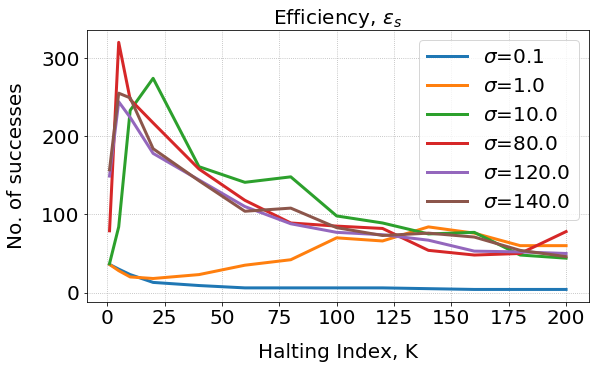}\label{fg:egg_ksenB}}
 
	\caption{Performance and efficiency results for the BH-S applied to the Egg-holder function for various combinations of $K$ and $\sigma$. A CPU time budget of 300s was applied to all simulations.  
	} 
	\label{fg:egg_ksen}
\end{figure} 

Figure~\ref{fg:egg_ksen} illustrates the impact on reliability and efficiency of the choice of halting index $K$. From Algorithm~\ref{bhmssalgo}, the maximum linear distance covered by the skipping procedure is $\sum_{k=1}^K R_k$, where each $R_k$ is distributed as the radial part of a centred Gaussian with standard deviation $\sigma$. This suggests that $K$ should not be too small, and the plot of efficiency in Figure~\ref{fg:egg_ksenB} indicates that $K$ should be at least 5 in our example (by default we take $K=25$).

It is seen that provided ($K \leq 5$), increasing $K$ tends to increase reliability while decreasing efficiency. This reflects the fact that larger $K$ allows the skipping procedure to travel further, thus increasing the likelihood of a direct transition to the global minimum basin, after which the BH-S algorithm would stop due to its monotonicity. In this way, greater $K$ increases reliability. On the other hand, greater $K$ increases the length of unsuccessful skipping trajectories. That is, each time the perturbed state $Y_n$ of Algorithm~\ref{bhmssalgo} is not accepted (after the local minimisation step of Algorithm~\ref{BHalgo}), the landscape is evaluated up to $K$ times without  advancing the optimisation. This implies that increased $K$ also typically leads to decreased efficiency. 

The tuning considerations discussed above for the BH-S algorithm can be summed up as follows. It should first be checked that $\sigma$ is large enough that the initial perturbation regularly falls outside the basin of the current state $X_n$. Having selected $\sigma$, $K$ should then be taken large enough that the skipping procedure regularly enters the sublevel set $C_n$. A practical suggestion here is to choose $K$ so that $K \sigma$ exceeds the diameter of the domain $D$.
 
\begin{figure}[!ht]
\centering
    \subfloat[][Percent of samples in desired basin: $\rho$]{\includegraphics[width=0.47\textwidth]{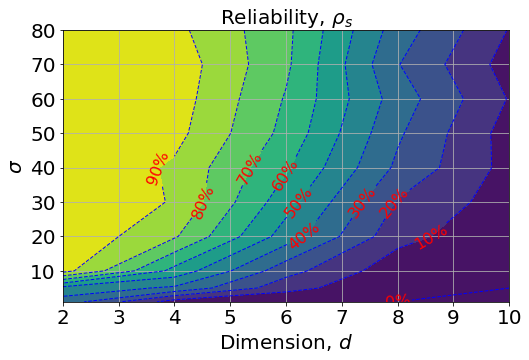} \label{fg:sch_pibA}} \quad
    \subfloat[][Total trials in desired basin: $\epsilon$]{\includegraphics[width=0.47\textwidth]{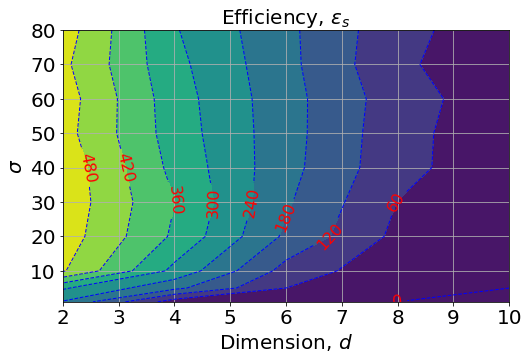} \label{fg:sch_pibB}} \quad

	\caption{Performance of the BH-S algorithm on the Schwefel-07 function for different combinations of domain dimension $d$ and perturbation variance $\sigma$. Note: the halting index was set to $K = 50$ with a CPU time budget of 300s for all simulations.}
	\label{fg:sch_pib}
\end{figure}

Figure~\ref{fg:sch_pib} confirms these guidelines in higher dimensions, by plotting the BH-S reliability and efficiency in dimensions up to 10 as $\sigma$ varies with the fixed choice $K=50$. It confirms that these performance metrics are relatively robust to the value of $\sigma$, provided that $\sigma$ is sufficiently large.
 
\subsection{Alternating BH-S and BH}
\label{sub:alt}
\begin{figure}[!ht]
 \centering
    \subfloat[][]{\includegraphics[width=0.46\textwidth]{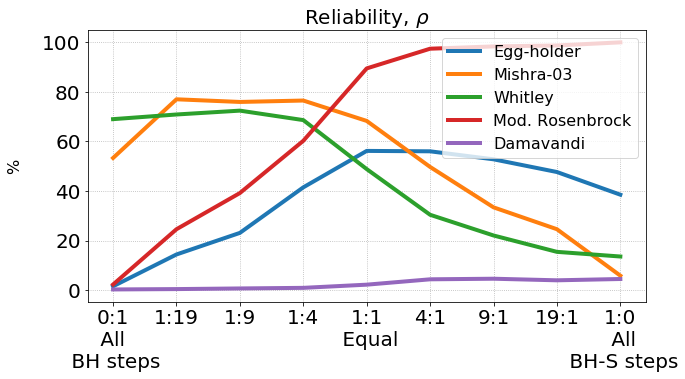} \label{fg:alt_pib}} \qquad
    \subfloat[][]{\includegraphics[width=0.46\textwidth]{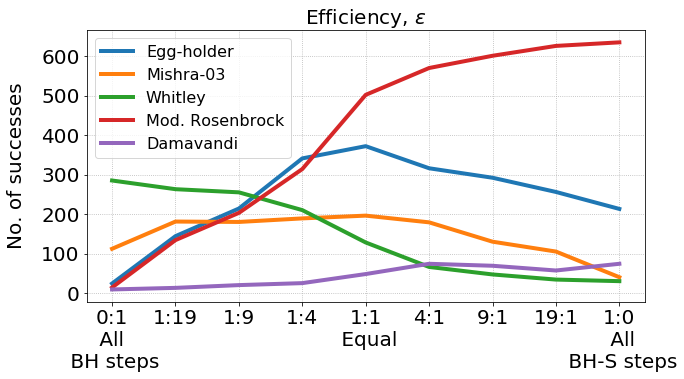} \label{fg:alt_pert}}

	\caption{Performance of the hybrid algorithm with varying proportions of BH to BH-S steps.}
 \label{fg:alt_graph}	
\end{figure}

In this section we explore a hybrid approach which is intended to overcome the challenges identified in Section~\ref{sub:nondistbas} for the monotonic BH-S algorithm by regularly including non-monotonic BH steps. Figure~\ref{fg:alt_graph} plots the reliability and efficiency metrics for this hybrid algorithm on various landscapes, as the ratio between BH-S and BH steps varies.

It can be seen that for the Mishra-03 and Whitley functions of Section~\ref{sub:nondistbas}, this hybrid improves both reliability and efficiency compared to BH-S. Indeed, the performance of a 1:1 ratio of BH and BH-S steps is comparable to that of BH for these landscapes. Further, on the landscapes of Section~\ref{sub:distbas}, this 1:1 ratio achieves performance superior to that of BH and somewhat comparable to that of BH-S. Thus if little is known about the problem's energy landscape \textit{a priori}, these results indicate that the 1:1 hybrid is to be preferred.
 
\section{Discussion and future work}
\label{sec:extens}

Basin hopping with skipping (BH-S) is a global optimisation algorithm inspired by both the basin hopping algorithm and the skipping sampler, an MCMC algorithm. As such, the MCMC literature also suggests potential extensions of this work. In adaptive MCMC parameter tuning is an online procedure driven by the progress of the chain \cite{yves2005}. A similar idea has been proposed for BH in~\cite{Gehrke2009} and is part of the of the SciPy implementation of the BH method. We believe it could be interesting as future work to devise an adaptive scheme for the halting index $K$ and the standard deviation $\sigma$, possibly reducing in this way the amount of tuning required to implement BH-S.
    
During the preparation of this paper we also explored the idea of sampling several directions and skipping in all of them simultaneously. As a negative finding, we report that preliminary results were clear that computational effort is best spent searching over a single, rather than multiple, directions. Our heuristic explanation is that the line is the shortest route between two sets, and so is the most efficient way to cover distance. An alternative, more sophisticated approach would be to introduce multiple BH-S particles which explore the energy landscape in a coordinated way. This could for instance be inspired by selection-resampling procedures as in sequential Monte Carlo sampling \cite{Doral2006}, or by an optimisation procedure such as particle swarm optimisation \cite{ken1995}.

\bibliography{libraryoptim}

\section*{Appendix}
\label{appendix}
Table~\ref{tb:results} records the results for all landscapes in Figure~\ref{fg:pibeps}. For each landscape, both BH and BH-S were hand tuned in order to maximise their efficiency as defined in Section~\ref{sub:method}. The following notation is used: 
\begin{itemize}
    \item $\rho_c$ and $\rho_s$ are the reliability of BH and BH-S respectively; 
    \item $\epsilon_c$ and $\epsilon_s$ are the efficiency of BH and BH-S respectively; 
    \item $\upsilon_1$ is the expected mean jump distance among random walk steps; 
    \item $\upsilon_s$ is the expected mean jump distance among skipping transitions, i.e., when $k>1$; 
    \item$\mathbb{P}_s$ is the probability that, conditional on the BH-S perturbation being accepted,  skipping had occurred ($k>1$);
    \item $\nu_1$ and $\nu_s$ are the expected mean jump distances among random walk steps ($k =1$) and skipping steps ($k>1$), respectively. 
\end{itemize}
 
\setlength\tabcolsep{6pt}	
\begin{longtable}[H]{ccccc||ccccccc}

\caption{Performance metrics for the BH and BH-S algorithms on all landscapes in Figure~\ref{fg:pibeps}.}\\
\label{tb:results}
\\ \hline\noalign{\smallskip}	

&\multicolumn{4}{c||}{BH}& \multicolumn{7}{c}{BH-S} \\ 
		\noalign{\smallskip}\hline\noalign{\smallskip}
		
{Function}&$\sigma$&$\rho_c$&$\epsilon_c$&$\nu_1$&$\sigma$&$K$&$\rho_s$&$\epsilon_s$&$\nu_1$&$\nu_s$&$\mathbb{P}_s$\\ \hline
\endfirsthead

\\ \hline\noalign{\smallskip}	

&\multicolumn{4}{c||}{BH}& \multicolumn{7}{c}{BH-S} \\ 
		\noalign{\smallskip}\hline\noalign{\smallskip}

{Function}&$\sigma$&$\rho_c$&$\epsilon_c$&$\nu_1$&$\sigma$&$K$&$\rho_s$&$\epsilon_s$&$\nu_1$&$\nu_s$&$\mathbb{P}_s$\\ \hline
\endhead
\noalign{\smallskip}\hline\noalign{\smallskip}

    \makecell{Carrom Table\\
    \begin{minipage}{.2\textwidth}
    \includegraphics[width=\linewidth]{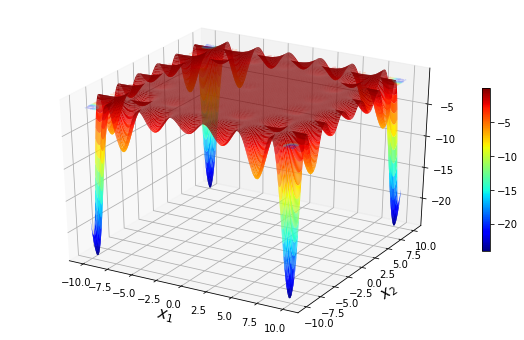}
    \end{minipage}}&	2	&	99.8	&	556	&	1.5	&	$\sqrt{2}$	&	10	&	100	&	925	&	2.9	&	9.4	&		86.8	\\ \hline
    
    \makecell{Cross in Tray\\
    \begin{minipage}{.2\textwidth}
    \includegraphics[width=\linewidth]{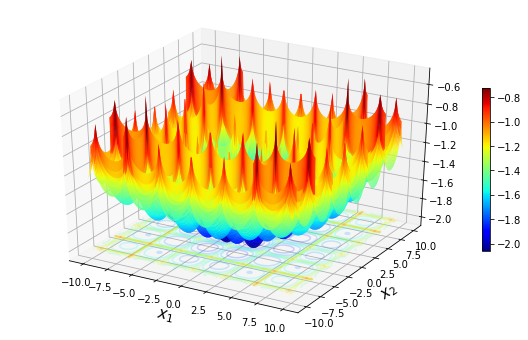}
    \end{minipage}}& 2	&	96.5	&	446	&	1.7	&	$\sqrt{2}$	&	10	&	100	&	676	&	2.7	&	9.3	&		85.1\\ \hline
    
    \makecell{Cross Leg Table\\
    \begin{minipage}{.2\textwidth}
    \includegraphics[width=\linewidth]{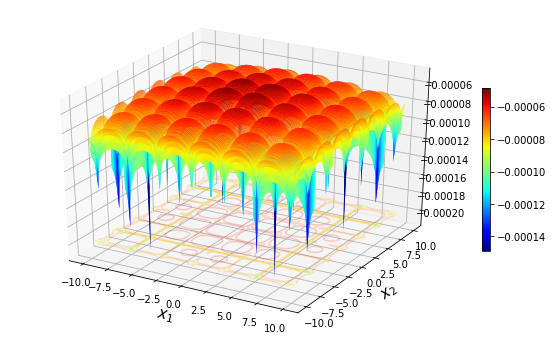}
    \end{minipage} }&0.4	&	15.5	&	48	&	0.8	&	$\sqrt{2}$	&	10	&	12.7	&	45	&	1.8	&	6.6	&	45.7			
\\ \hline
    
    \makecell{Damavandi\\
    \begin{minipage}{.2\textwidth}
    \includegraphics[width=\linewidth]{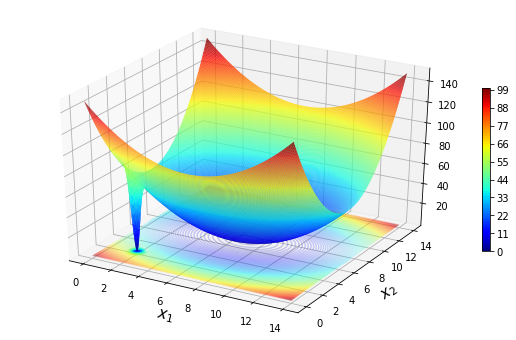}
    \end{minipage}}&0.1	&	0.2	&	3	&	0.5	&	0.3	&	150	&	32.9	&	28	&	N/A	&	34.9	&	100			
 \\ \hline
    
    \makecell{Eggcrate\\
    \begin{minipage}{.2\textwidth}
    \includegraphics[width=\linewidth]{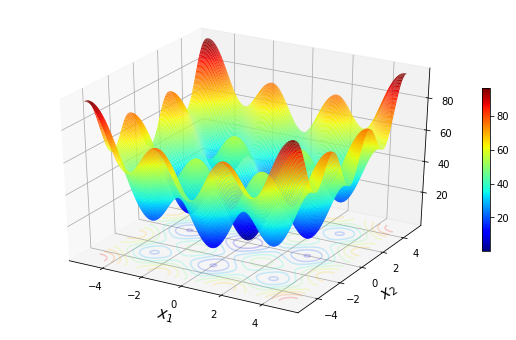}
    \end{minipage}}&	1	&	99.7	&	377	&	1	&	1	&	10	&	99.7	&	647	&	1.9	&	7.6	&	94.8			
\\ \hline
    
    \makecell{Egg-holder\\
    \begin{minipage}{.2\textwidth}
    \includegraphics[width=\linewidth]{eggholder.png}
    \end{minipage}}&100	&	2.2	&	13	&	12.5	&	10	&	25	&	38.7	&	116	&	7.1	&	178.1	&	98.6
    \\ \hline
    
    \makecell{El Attar Vidyasagar Dutta\\EAVD\\ 
    \begin{minipage}{.2\textwidth}
    \includegraphics[width=\linewidth]{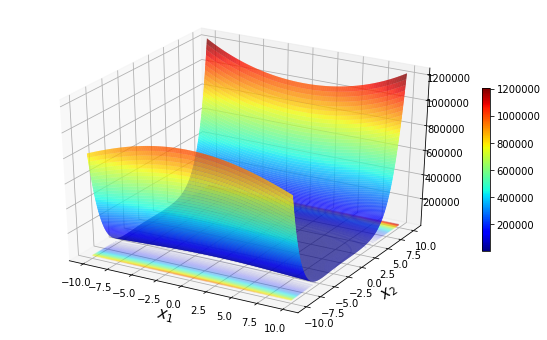}
    \end{minipage}}& 8	&	99.6	&	231	&	3.2	&	5	&	10	&	63.4	&	393	&	4.8	&	5.8	&	39.3		\\ \hline
    
    \makecell{Freudenstein-Roth\\
    \begin{minipage}{.2\textwidth}
    \includegraphics[width=\linewidth]{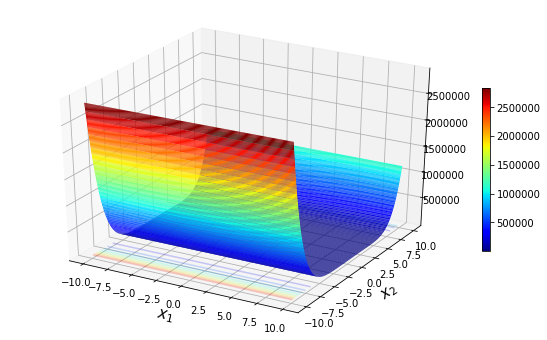}
    \end{minipage}} &2	&	82.6	&	176	&	1.6	&	$\sqrt{2}$	&	10	&	97.2	&	551	&	4.9	&	11.3	&	99.5			
 \\ \hline
    
    \makecell{Holder Table\\
    \begin{minipage}{.2\textwidth}
    \includegraphics[width=\linewidth]{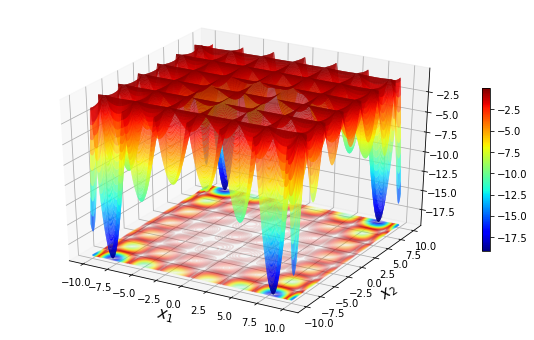}
    \end{minipage}}& 2	&	99.8	&	453	&	1.4	&	$\sqrt{2}$	&	10	&	100	&	695	&	2.4	&	9.1	& 82.4			\\ \hline
    
    \makecell{Keane\\
    \begin{minipage}{.2\textwidth}
    \includegraphics[width=\linewidth]{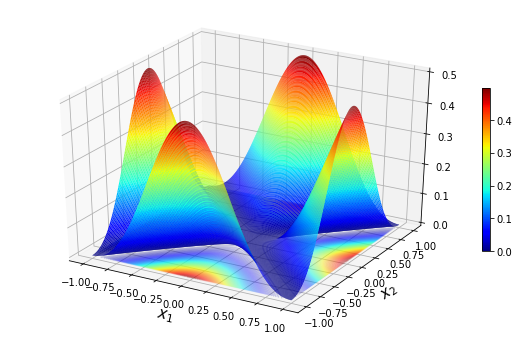}
    \end{minipage}}& 2	&	47	&	272	&	1.8	&	0.9	&	25	&	53.6	&	301	&	1.6	&	12.4	&	96.6	\\ \hline
  
     \makecell{Mishra-03\\
    \begin{minipage}{.2\textwidth}
    \includegraphics[width=\linewidth]{mishra03.png}
    \end{minipage} }& 2	&	65.9	&	56	&	1.8	&	$\sqrt{2}$	&	10	&	5.4	&	17	&	1.6	&	12.1	&		96.7 \\\hline
    
    \makecell{Modified Rosenbrock\\
    \begin{minipage}{.2\textwidth}
    \includegraphics[width=\linewidth]{modrosenbrock.png}
    \end{minipage} }&0.4	&	5.3	&	1	&	0.8	&	0.4	&	25	&	83.8	&	31	&	N/A	&	7.7	&		100	\\\hline

    \makecell{Price 02\\
    \begin{minipage}{.2\textwidth}
    \includegraphics[width=\linewidth]{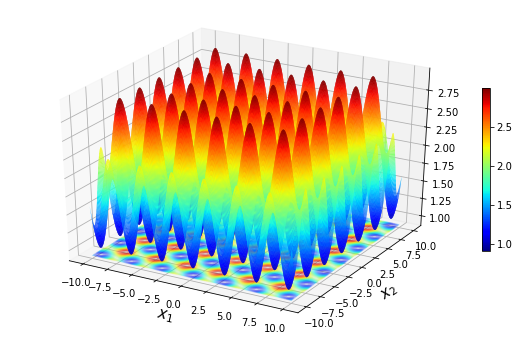}
    \end{minipage}}&2	&	44.4	&	234	&	1.8	&	0.9	&	25	&	60.6	&	208	&	N/A	&	17.1	&	100			
\\ \hline
    
    \makecell{Rana\\
    \begin{minipage}{.2\textwidth}
    \includegraphics[width=\linewidth]{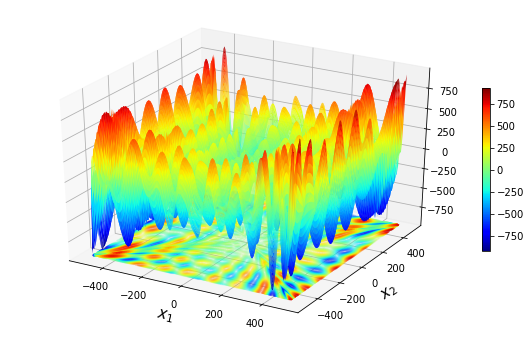}
    \end{minipage}}& 200	&	5.5	&	13	&	17.6	&	5	&	75	&	20.5	&	18	&	1.8	&	224.5	&	98.5			
\\\hline

    \makecell{Rosenbrock\\
    \begin{minipage}{.2\textwidth}
    \includegraphics[width=\linewidth]{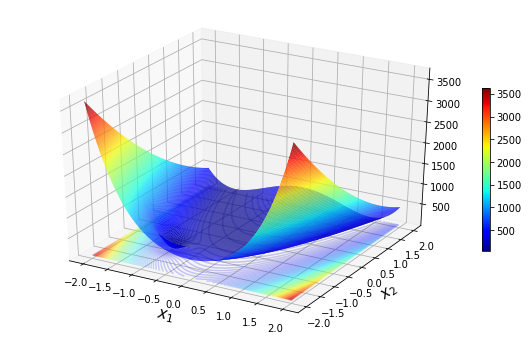}
    \end{minipage} }&	0.2	&	99.3	&	149	&	0.6	&	2	&	10	&	100	&	695	&	N/A	&	N/A	&	0			
\\ \hline
    
    \makecell{Schwefel-07\\
    \begin{minipage}{.2\textwidth}
    \includegraphics[width=\linewidth]{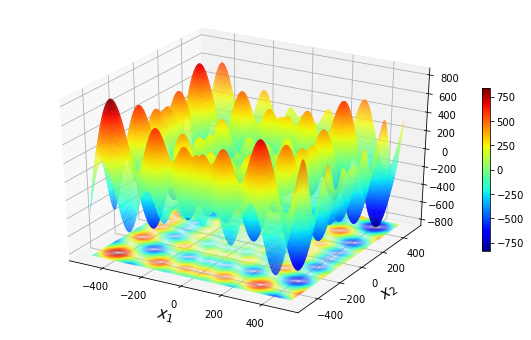}
    \end{minipage}}& 10	&	4.1	&	38	&	4	&	7	&	25	&	61.9	&	275	&	N/A	&	181.6	&	100			
\\ \hline
    
    \makecell{Styblinski Tang\\
    \begin{minipage}{.2\textwidth}
    \includegraphics[width=\linewidth]{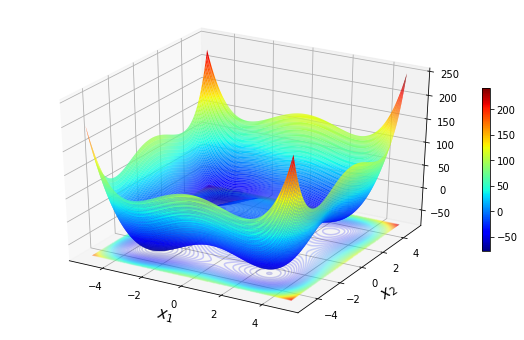}
    \end{minipage}}& 1	&	99.6	&	537	&	1.2	&	1	&	10	&	100	&	840	&	N/A	&	8.7	&	99.8			
	\\ \hline
  
    \makecell{Whitley\\
    \begin{minipage}{.2\textwidth}
    \includegraphics[width=\linewidth]{whitley.png}
    \end{minipage}}& 0.4	&	86.4	&	121	&	0.7	&	0.7	&	50	&	31.8	&	21	&	0.6	&	17	&	37.7			
\\ \hline
    
    \makecell{Zirilli\\
    \begin{minipage}{.2\textwidth}
    \includegraphics[width=\linewidth]{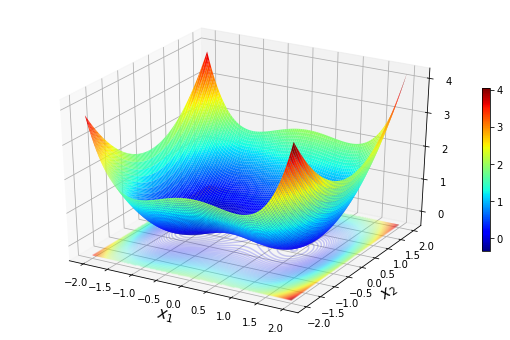}
    \end{minipage}}& 0.2	&	99.9	&	707	&	0.6	&	$\sqrt{2}$ 	&	10	&	97.9	&	987	&	1.9	&	5.3	&	49.4			
		\\ \hline

\end{longtable}

\end{document}